\documentclass[letterpaper, 10 pt, conference]{ieeeconf}  

\IEEEoverridecommandlockouts                              
\usepackage{graphics} 
\usepackage{epsfig} 
\usepackage{amsmath} 
\usepackage{amssymb}  

\usepackage{verbatim}
\usepackage{graphicx} 
\usepackage{epsfig} 
\usepackage{epstopdf}
\usepackage{amsmath} 
\usepackage{amssymb}  
\usepackage{caption}
\usepackage{color}

\makeatletter
\newcommand{\figcaption}{\def\@captype{figure}\caption}
\newcommand{\tabcaption}{\def\@captype{table}\caption}
\makeatother

\title{\LARGE \bf
On Kelly Betting: Some Limitations
}

 \author{Chung-Han Hsieh$^{1}$ and B. Ross Barmish$^{2}$
\thanks{\hskip -10pt ${}^1$Chung-Han Hsieh is a graduate student working towards to his Ph.D. degree in the Department of Electrical and Computer Engineering, University of Wisconsin, Madison, WI 53706. E-mail: hsieh23@wisc.edu.}
\thanks{\hskip -10pt ${}^2$B. Ross Barmish is a faculty member in  the Department of Electrical and Computer Engineering, University of Wisconsin, Madison, WI 53706. \mbox{E-mail}: barmish@engr.wisc.edu.}}

\begin{document}

\maketitle
\thispagestyle{empty}
\pagestyle{empty}
\begin{abstract}
The focal point of this paper is the so-called {\it Kelly Criterion}, a prescription for optimal resource allocation 
among a set of gambles which are repeated over time. The criterion calls for maximization of the expected value of the
logarithmic growth of wealth. Considerable literature exists providing the rationale for such an optimization. This paper begins
by describing some of the limitations of the Kelly-based theory in the existing literature. 
To this end, we fill a void in published results by providing specific examples quantifying what can go wrong when Taylor-style approximations are used and when wealth drawdowns are considered. For the case of drawdown, we describe some research directions which we feel are promising for improvement of the theory. \\
\end{abstract}


\maketitle

\section {Introduction}
\noindent
The focal point of this paper is the so-called Kelly Criterion
introduced in the seminal paper~\cite{Kelly_1956}. Given $n$ gambles with return governed by some random vector~\mbox{$X \in \mathbb{R}^n$}, Kelly's theory indicates what fraction~$K_i$ of one's account value~$V$ to invest in the $i$-th bet. Letting~$K$ be the column vector with components~$K_i$, the classical formulation of this problem requires~$K_i \geq 0$ for~$i = 1,2,...,n$ and~$$
K_1 + K_2 + \cdots K_n \leq 1.
$$ 
The problem formulation also includes the standing assumption that this gamble  is repeated over and over again via independent and identically distributed (i.i.d.) trials for~$X$ and that~$K$ is such that survival is assured. This notion will be made precise in the sequel. 
\\
\\
Noting that the account value begins at some initial level~$V(0) > 0$ and letting~$X(k)$ be the~$k$-th outcome for~$X$, evolution to terminal state~$V(N)$ is described sequentially by the recursion
$$
V(k+1) = (1+K^TX(k))V(k).
$$
Letting~${\cal X} \subset \mathbb{R}^{n}$ denote the support of~$X$ which we  assume to be closed, in order to assure satisfaction of the survival requirement, admissible~$K$ must satisfy the condition
$$
\min_{X \in {\cal X}}K^TX \geq -1.
$$
Henceforth, to denote the totality of the constraints above, we write~$K \in {\cal K}$ and note that $\cal K$ is convex. To conclude this overview, it is noted that there are many possible variations and extensions of this problem formulation in the literature. For example, one can allow~$K_i > 1$ to include leverage considerations and~$K_i < 0$ to model short sales. Finally, we mention one of the most important application areas for the ideas to follow: trading and portfolio balancing problems in financial markets. Following the results in~\cite{Kelly_1956}, we see a trail in the literature over the subsequent decades dealing with all sorts of
applications, generalizations and improvements of the theory; e.g., see~\cite{Rotando_Thorp_1992}, \cite{Thorp_2006}, \cite{Maclean_Sanegre_Zhao_Ziemba_2004}, and~\cite{Nekrasov_2014}.
\\
\parindent = 0pt
\subsection{Problem Formulation}
The classical Kelly problem is to select~$K \in {\cal K}$ so as maximize the expected value of the logarithmic {\it growth}
$$
g(K) \doteq \frac{1}{N} \mathbb{E} \left[ {\log \left( {\frac{{V(N)}}{{V(0)}}} \right)} \right].
$$
Using the recursion for~$V(k)$ above, the additivity of the log function, the fact that $X(k)$ are i.i.d., it is easy to show that the expected log-growth function reduces to
\begin{align*}
	g(K) &= \mathbb{E}[\;\log(1 + K^T X)\;] \\
	&= {\int_{\mathcal{X}} {\log } (1 + K^T x){f_X}(x)dx}
\end{align*}
where~$f_X(x)$ denotes the probability density function for~$X$.
Subsequently, when the constraint~$K \in {\cal K}$ is included, it is easy to show that the optimal logarithmic growth
$$
g^* \doteq \max_{K \in {\cal K}}g(K)
$$
is a concave program in $K$.
\\
\\
 To provide one of the simplest possible illustrations for all of the above, the literature in~\cite{Kelly_1956} considers flipping a biased coin with gambling return~$X(k) = 1$ with probability~$p > 1/2$ and~$X(k) = -1$ with probability~$1-p$. In this scenario,~$f_X(x)$ is described by a pair of Dirac Delta functions and it is readily shown by straightforward differentiation of~$g(K)$ above that the optimal fraction,~\mbox{$K = K^*$}, is given by~$K^* = 2p-1$.
\\

\subsection{Why Use the Logarithmic Growth?}
Use of the Kelly Criterion has a number of advantages over the use of  the more classical expected value of terminal wealth~$\mathbb{E}[V(N)]$. To illustrate why this is so, for~$n=1$, if~$\mathbb{E}[X(k)]$  is just ``slightly'' positive, it is easy to see that the optimum is obtained by making~$K$  as large as permitted; e.g., for the case of an even-money bet on a  biased coin with winning probability~$p = 0.5 + \varepsilon$, no matter how small the advantage~$\varepsilon > 0$ is,  maximizing~$\mathbb{E}[V(N)]$ dictates using~$K=1$. Such as strategy is arguably far too aggressive to use for a game which is being played over and over again. With~$N$ large, it is almost certain that~$V(k)$ will be drawn down to zero; i.e., gambler's ruin will occur.\\
\\
In contrast to the use of~$\mathbb{E}[V(N)]$ above, the Kelly Criterion, in its use of~$\mathbb{E}[\log(V(N)]$, automatically factors some degree of risk into the analysis. For the case of the coin above with small~$\varepsilon > 0$, the optimum turns out to be~$K = 2\varepsilon$, thereby much more likely to avoid gambler's ruin. By taking into account the exponential growth rate of wealth and carrying out the myopic period-by-period optimization leading to optimal logarithmic growth, a number of desirable properties result thereby making the Kelly Criterion a  powerful tool in finance; see~\cite{Maclean_Thorp_Ziemba_2010} where a nice summary of both the desirable and undesirable properties are given.
In this regard, of foremost importance is the following: When the optimal Kelly fractions~$K_i$ increase, various risk measures become unacceptably large. Hence, the literature also includes a number of papers dealing with ``fractional strategies.'' Essentially, this amounts to reduction of the~$K_i$, often in ad hoc manner; e.g., see~\cite{Maclean_Sanegre_Zhao_Ziemba_1999}. Finally, to provide further context for the sections to follow, we mention other related papers in the literature, see~\cite{Hakansson_1972}-\cite{Maslov_Zhang_1998},~\cite{Maclean_Sanegre_Zhao_Ziemba_2004}, \cite{Rising_Wyner_2012}-\cite{Nekrasov_2014_book}, and single out~\cite{Calafiore_Monastero_2010} which has the same control-theoretic point of view described below.\\ 

\parindent = 0pt
\subsection{Feedback Control System Point of View}
\begin{figure}[htbp]
   \centering
   	\graphicspath{{fig/}}	\setlength{\abovecaptionskip}{0.1 pt plus 0pt minus 0pt}	   
   	\includegraphics[width=0.45\textwidth]{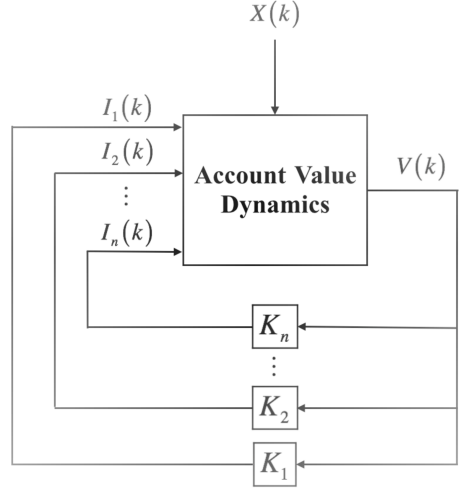}
   	\figcaption{Feedback Control Equivalent of Kelly Betting\\}
   	\label{fig.Block_diagram}
   \end{figure}

The problem formulation above is readily interpreted in terms classical feedback control theory. That is, we view~$V(k)$ as the state of a system with linear feedback and~$n$ inputs corresponding to the investment levels~$I_i(k)$ for each of the gambles. That is, the~$i$-th input of the control signal is given by
$$
I_i(k) = K_iV(k)
$$	
with~$K_i \geq 0$ viewed as a feedback gain. Subsequently, the state for this stochastic system is updated via the equation
\begin{align*}	
V(k+1) & = V(k) + \sum_{i = 1}^{n}I_i(k)X_i(k)V(k)\\
		& = (1 + K^TX(k)) V(k).
\end{align*}
This type of feedback-control configuration is depicted in Figure~\ref{fig.Block_diagram}; see~\cite{Barmish_Primbs_2015} where this paradigm is pursued in much greater detail.\\

\subsection{Plan for Sections to Follow}
Although limitations of the Kelly-based theory are mentioned in the existing literature, there is a paucity of specific examples illustrating the degree to which things can ``go wrong.'' To this end, Section~2 concentrates on approximation methods in the literature which are used to optimize the allocation vector~$K$. As shown in the sequel, the Taylor series are used to approximate the log-growth function, we see that the solution which is obtained may be either infeasible or lead to performance which is significantly lower than that of the true optimum. In addition, we show that approximate solutions may have a certain ``inefficiency property'' which is undesirable.
\\
\\
Although our examples to follow provide specific realizations of the ``badness'' which can occur, it should also be noted that a ``remedy'' is readily available. That is, some papers, for example, see \cite{Cover_1984} and \cite{Maclean_Sanegre_Zhao_Ziemba_2004}, recognize that the log-growth problem is a concave program. 
Hence, it is arguable that approximation methods are not needed because there are readily available commercial codes which efficiently solve the problem at hand; see~\cite{Boyd_2004} and~\cite{Grant_Boyd_2015}. At the time that some of the earlier papers were written, such codes were not readily available and authors either resorted to approximation or developed algorithms of their own; see~\cite{Cover_2012}.\\
\\
In Section 3, a concern which is much more serious than approximation is addressed --- the issue of wealth drawdown. Suffice it to say, the literature already recognizes that the Kelly gains~$K_i$ which result, although being log-growth optimal, may be too aggressive in the short term;~i.e., the wealth level~$V(k)$ may fall to unacceptably low levels along sample pathes. For this reason, as mentioned earlier, some authors resort to a so-called ``fractional'' betting scheme by scaling back the~$K_i$; e.g., see~\cite{Maclean_Sanegre_Zhao_Ziemba_1999} and~\cite{Maclean_Thorp_Ziemba_2010}. Other authors resort to incorporation of constraints to reduce the drawdown effect; e.g., see~\cite{Maclean_Sanegre_Zhao_Ziemba_2004}. After quantifying some of the negatives regarding drawdown, in Section 4, we describe some research directions which we feel are promising for mitigation of the drawdown problem. Finally, in Section 5, some conclusions are given and other directions of research are mentioned.  \\

%
\parindent = 0pt
\section{Negatives Associated with Approximation}
In order to obtain the optimal logarithmic growth rate~$g^*$ above, as previously mentioned, one approach in the literature involves approximation --- either a multivariate Taylor expansion to the log-growth function is used or~$X(k)$ is treated as a Geometric Brownian Motion and low order expansion terms are used; e.g., see~\cite{Thorp_2006}, \cite{Rising_Wyner_2012}, \cite{Nekrasov_2014}, and~\cite{Nekrasov_2014_book}. The main objective in this section is to point out some ``pitfalls'' associated with approximate solution. While approximation-based closed-form solutions for the optimal~$K$ provide a degree of insight into the risk-return tradeoffs, concrete examples do not appear in the literature which demonstrate scenarios where approximation methods lead to erroneous results. Suffice it to say, when the range of variation of~$X(k)$ can be large, the true optimum~$K = K^*$ and associated logarithmic growth~$g(K^*)$ can differ considerably from its approximation.
\\
\subsection{Example Involving Approximation}
We consider the somewhat attractive gamble for which~\mbox{$n = 1$},~$X = 0.15$ with probability~$p = 0.95$ and~$X = -0.95$ with probability~$p = 0.05$. We call this bet ``attractive" in a central-limiting sense; i.e., since~$\mathbb{E}[X] = 0.095$, repeated i.i.d. trials, will almost certainly lead to success. Now, according to~\cite{Nekrasov_2014}, using the approximation
$$
\mathbb{E}[\log(1 + KX)] \approx K \mathbb{E}[X] - \frac{1}{2}K^2 \mathbb{E}[X^2],
$$
it is straightforward to see that the associated optimum investment fraction $K$, call it $\scalebox{1.5}{$\kappa$}_{\rm Taylor}$, is given by
\[
 {\scalebox{1.5}{$\kappa$}}_{\rm Taylor}= \frac{\mathbb{E}[X]}{\mathbb{E}[X^2]} = 1.4286.
\]
Note that this solution is not feasible because $K \in [0,1]$ is assumed. Hence, to guarantee feasibility, a saturation function is introduced for the approximate solution above. Thus, the optimal approximate solution with saturation, call it $K_{\rm Taylor}$, and the associated expected log-growth are as follows   
\begin{align*}
& K_{\rm Taylor} = SAT \left[ \scalebox{1.5}{$\kappa$}_{\rm Taylor} \right] = 1;\\
& g(K_{\rm Taylor}) \approx -0.017
\end{align*}
where $SAT[x]$ is a saturation function; i.e., for $x <0$, $SAT[x] = 0$; for $0 \le x \le 1$, $SAT[x]=x$ and for $x>1$, we have $SAT[x]=1$. 
\\
\\
An alternative approach, for example, see~\cite{Thorp_2006} and \cite{Rising_Wyner_2012}, with~$X(k)$ being treated as a Geometric Brownian Motion with drift~$ \mu = \mathbb{E}[X]$ and variance~$\sigma^2 = VAR(X)$. A subsequent Taylor approximation leads to approximate solution, call it~$\scalebox{1.5}{$\kappa$}_{\rm GBM}$, as
	\[
	\scalebox{1.5}{$\kappa$}_{\rm GBM} = \frac{\mathbb{E}[X]}{VAR[X]} = 1.6529.
	\]
Similarly, it is infeasible with restriction  $K \in [0,1]$ so the saturation is required. Here, the associated optimal approximate solution, call it as $K_{\rm GBM}$, and the corresponding expected log-growth are given by  
\begin{align*}
	& K_{\rm GBM} = SAT \left[ \scalebox{1.5}{$\kappa$}_{\rm GBM}  \right] =1; \\
	& g(K_{\rm GBM}) \approx -0.017.
 \end{align*}
In contrast to the two approximate solutions above, the true optimum, as described in Section~1, is obtained by maximizing the expected logarithmic growth
	$$
	g(K) =  0.95\log(1 + 0.15K)    + 0.05 \log( 1 - 0.95K)
	$$
which, by straightforward differentiation, leads to a feasible solution in~$[0,1]$ and optimal growth given by
	$$
	K^* \approx 0.6667;\;\; g(K^*) \approx 0.0404.
	$$
A summary of all three solutions is given in Figure~2. Ironically, the approximation-based results yields the minimum growth of $g(K)$ rather than the desired maximum. Suffice it to say, the combination of approximation and saturation due to constraint violation can lead to significant error.
	
	\begin{figure}[htbp]
		\centering
		\graphicspath{{fig/}}	\setlength{\abovecaptionskip}{0.1 pt plus 0pt minus 0pt}	   
		\includegraphics[width=0.45\textwidth]{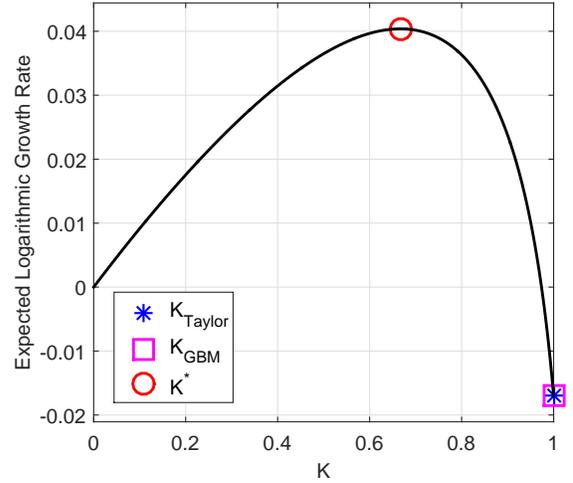}
		\figcaption{Expected Logarithmic Growth Rate \\}
		\label{fig:GvsK}
	\end{figure}
	
To quantify further, we convert the expected log-growth into an annualized rate of return using
	\[
	r(K) \doteq \frac{1}{\Delta t} (e^{g(K)} - 1).
	\]
Assuming daily betting, where $\Delta t$ is the time between bets in years, we take $\Delta t = 1/252$ and then, the corresponding expected annualized rates of return are computed to be 
\begin{align*}
	& r(K^*) \approx 10.384;\\
	& r(K_{\rm Taylor}) = r(K_{\rm GBM}) \approx -3.443.
\end{align*}
In other words, the approximate betting schemes perform poorly compared to what is possible.
\\


\parindent = 0pt
\subsection{More Realistic Example with Real Stock Data}
\begin{figure}[htbp]
	\centering
	\graphicspath{{fig/}}	\setlength{\abovecaptionskip}{0.1 pt plus 0pt minus 0pt}	   
	\includegraphics[width=0.45\textwidth]{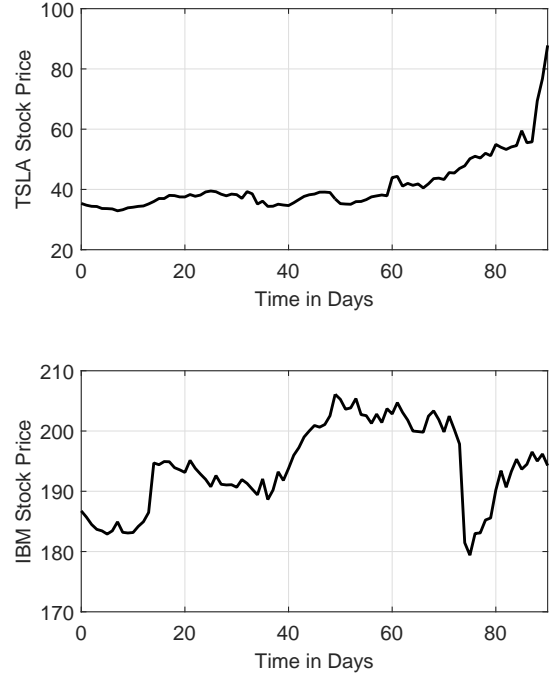}
	\figcaption{Two Stock Prices: TSLA and IBM\\}
	\label{fig3:Two_Stock_Prices}
\end{figure}

In this example, we further consider the problems associated with approximation by using an example involving real data for two stocks: Tesla Motors and IBM during the ninety-day period January 2, 2013 until May 13, 2013. We used the adjusted daily closing prices, see Figure~\ref{fig3:Two_Stock_Prices}, to estimate the joint probability mass function and carried out an in-sample constrained maximization of~$g(K)$ subject to the constraint~$K_1 + K_2 \le 1$, we obtain
$$
K_1^* = 1;\; K_2^* = 0.
$$
That is, the optimum log-growth solution involves all funds invested in Tesla and no investment in IBM. Now, suppose instead that one computes the Taylor-based solutions; i.e.,
\begin{align*}
	& \scalebox{1.5}{$\kappa$}_{\rm Taylor} = \Sigma^{-1} (X) \mathbb{E}[X] = [5.321\;\; 2.725]^T;\\
	& \scalebox{1.5}{$\kappa$}_{\rm GBM} = \bar \Sigma^{-1} (X) \mathbb{E}[X] = [5.599\;\; 2.681]^T
\end{align*}
where $\Sigma(X)$ is the second moment matrix for $X$ and $\bar \Sigma(X)$ is the covariance matrix for $X$.
Note that the approximate solutions $\scalebox{1.5}{$\kappa$}_{\rm Taylor}$ and $\scalebox{1.5}{$\kappa$}_{\rm GBM}$ are infeasible since the constraint is violated. The true optimum solution and approximate solutions along with $\scalebox{1.5}{$\kappa$}_{\rm Taylor}$ and $\scalebox{1.5}{$\kappa$}_{\rm GBM}$ are seen in Figure~\ref{fig:GvsK}. 
\\
\\
Given the constraint violation $K_1 + K_2 >1$ for the Taylor and GBM solutions, one standard approach is to project these solutions onto the constraint satisfaction set. That is, we take
\begin{align*}
	& K_{\rm Taylor} = Proj ( \scalebox{1.5}{$\kappa$}_{\rm Taylor}) \approx [0.661\;\; 0.339]^T; \\
	& K_{\rm GBM}  = Proj (\scalebox{1.5}{$\kappa$}_{\rm GBM}) \approx [0.661 \;\; 0.339]^T
\end{align*}
where $Proj(\cdot)$ is a projection function given by
$$
Proj(K_1, K_2) \doteq \left[\frac{K_1}{K_1+K_2} \;\; \frac{K_2}{K_1+K_2} \right ]^T
$$ for all nonnegative $K_1,K_2$ with not both  $K_1,K_2 = 0$. However, one should note that although the projection procedure provides a way for yielding a feasible solution, the projected solution may not be the optimal. \\

\begin{figure}[htbp]
	\centering
	\graphicspath{{fig/}}	\setlength{\abovecaptionskip}{0.1 pt plus 0pt minus 0pt}	   
	\includegraphics[width=0.45\textwidth]{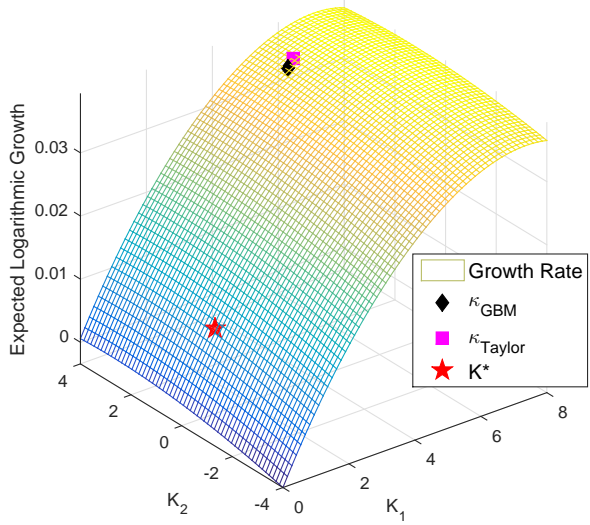}
	\figcaption{Constraint Violation Example for Two Stocks Case\\}
	\label{fig:GvsK}
\end{figure}

\parindent=0pt

\subsection{Inefficiency of Approximate Solution}
\begin{figure}[htbp]
	\centering
	\graphicspath{{fig/}}	\setlength{\abovecaptionskip}{0.1 pt plus 0pt minus 0pt}	   
	\includegraphics[width=0.45\textwidth]{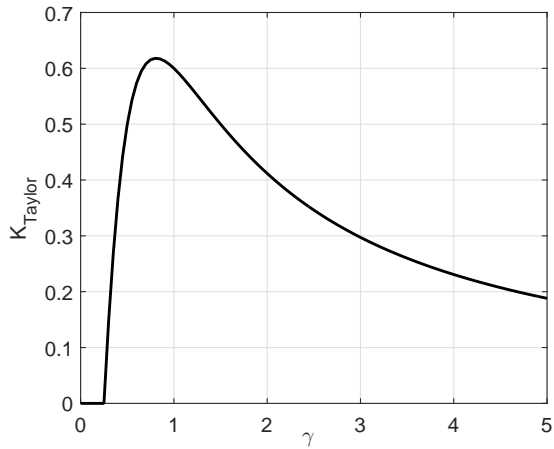}
	\figcaption{$K_{\rm Taylor}(\gamma)$ Plot for~$ p = 0.8$\\}
	\label{fig: Irrationality}
\end{figure}

In this subsection, we point out another danger associated with the use of approximate solutions. 
The takeoff point is the following principle widely used in finance: If two investments have the same risk, the one with the smaller reward will be discarded and deemed to be ``inefficient.''
We claim that the use of either the approximation~$K_{\rm Taylor}$  or~$K_{\rm GBM}$ might be inefficient.
We now provide such an example using~$K_{\rm Taylor}$ and note that the same example can be used for~$K_{\rm GBM}$ too. Indeed,  we consider a random variable~$X$ described as follows: Given~$\gamma > 0$, we have~\mbox{$X = \gamma$} with probability~$p>0$ and~$X = -1$ with probability~$1-p$.
Using the Taylor approximation, as a function of reward level~$\gamma$, a straightforward calculation yields
$$
K_{\rm Taylor}(\gamma) = SAT\left[\frac{p\gamma + p-1}{p\gamma^2  - p +1}\right].
$$
In order for~$K_{\rm Taylor}(\gamma)$ to be efficient from an economic risk-taking point of view, it should have the following property: When~$\gamma_2 \geq \gamma_1 \geq 0$,  we require $K(\gamma_2) \geq K(\gamma_1).$ That is, if the bet associated with~$\gamma_2$ offers more reward with the same probabilities of success and failure as those for~$\gamma_1$, a rational gambler should invest at least as much in the~$\gamma_2$ bet as the~$\gamma_1$ bet.  We claim that the Taylor-based approximation scheme fails to satisfy this condition. 
To establish  this  claim, it suffices to show that~$dK_{\rm Taylor}/d\gamma$ can be negative with the $K_{\rm Taylor}(\gamma)$ in~$(0,1)$. Indeed, we calculate
	$$
	\frac{dK_{\rm Taylor}}{d\gamma} = -p\;\frac{p\gamma^2 + 2(p-1)\gamma + p-1}{{\left(p\, \gamma^2 - p + 1\right)}^2}
	$$
and note that the denominator cannot vanish. Hence, we see~$dK_{\rm Taylor}/d\gamma < 0$ for
	$$
	\gamma > \gamma^*(p) \doteq \frac{1-p + \sqrt{1-p}}{p}
	$$
which corresponds to the zero-crossing of the numerator. In the Figure~\ref{fig: Irrationality}, the plot of $K_{\rm Taylor}$ is given for $p = 0.8$. It is readily apparent that the claimed inefficiency occurs for parameter range $\gamma > \gamma^*(0.8) \approx 0.809. $ \\

\section{Negatives Associated with Drawdown}
As described in \cite{Maclean_Sanegre_Zhao_Ziemba_2004} and \cite{Malekpour_Barmish_2013}, control of drawdown, that is,  control of the drops in wealth from peaks to subsequent lows, is one of great concern from a risk management perspective. 
In this section, we first demonstrate that Kelly betting often results in very poor drawdown performance. Then, we discuss some approaches for mitigating the drawdown problem within this framework. Along any sample path~$V(k)$, the {\it maximum percentage drawdown} is defined as
\[
D(K) \doteq \max_{0 \le l \le k \le N} \frac{V(l) - V(k)}{V(l)}.
\]
In the sequel, we often drop the word ``percentage" for expression simplicity. 
\\
\\
As mentioned in Section 1, use of Kelly fraction may lead to a significant drawdown since it is too aggressive. To quantify see how bad the drawdown can be, consider betting~$N$ times of single coin flipping gamble for which~$X = 1$ with probability~$p$ and~$X=-1$ with probability $1-p$, then it is easy to show that the probability of maximum drawdown greater than or equal to any fraction~$K \in (0,1)$ is given by
$$
P(D(K) \ge K) = 1 - p ^N.
$$
Now if we take $N=252$ and $p= 0.99$, using the optimal Kelly investment fraction $K^* = 2 p -1 = 0.98,$ It follows that there is a $ 92 \% $ chance that maximum drawdown is over $98\% $. That is, there is a large drawdown occurs with very high probability. A similar analysis using the Markov inequality also leads to the same conclusion. 
\\
\subsection{Control of Drawdown}
 To control the drawdown, one possible choice is to add probabilistic constraint to the optimization of log-growth;~i.e., given $0< \varepsilon <1$ and $0< \delta <1$, consider the constraint 
\[
P(D(K) \le \varepsilon) \ge 1 - \delta.
\]
Alternatively, instead of using the probabilistic  constraint above, we can use the expected maximum drawdown; i.e., given $0<\varepsilon < 1$, consider the drawdown constraint as
\[
\mathbb{E}[D(K)] \le \varepsilon.
\] 

We now revisit the example used in Section~2 with~\mbox{$n=1$}, \mbox{$X = 0.15$} with probability~$p=0.95$ and~$X=-0.95$ with probability~$p=0.05$. 

\begin{figure}[htbp]
	\centering
	\graphicspath{{fig/}}	\setlength{\abovecaptionskip}{0.1 pt plus 0pt minus 0pt}	   
	\includegraphics[width=0.45\textwidth]{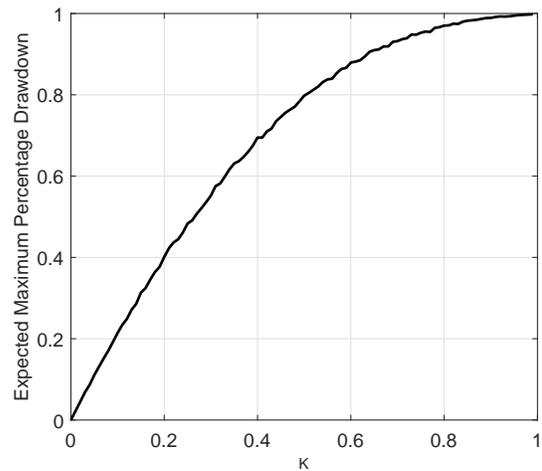}
	\figcaption{Expected Maximum Percentage Drawdown Versus $K$\\}
	\label{fig5:E_DD_vs_K_for_motivation_ex}
\end{figure}

Using the optimum fraction~$K^* = 0.6667$ already found, it is clear to see from Figure~\ref{fig5:E_DD_vs_K_for_motivation_ex} that the corresponding expected maximum drawdown are
	\begin{align*}
	& \mathbb{E}[D(K^*)]=\mathbb{E}[D(0.6667)] \approx 0.903 \;; \\
	& \mathbb{E}[D(K_{\rm Taylor})] = \mathbb{E}[D(K_{\rm GBM})] = \mathbb{E}[D(1)]\approx 1.0 \; .
	\end{align*}
This shows that the approximation solution leads to an almost sure ruin. Now, suppose the gambler adds constraint~$\mathbb{E}[D(K)] \le 0.2$. Then, based on Figure~\ref{fig5:E_DD_vs_K_for_motivation_ex}, it is clear to see that the optimal investment fraction $K$ reduces to~$K = K^* \approx 0.1.$
\\
\section{Research Directions Involving Drawdown} 
\begin{figure}[htbp]
	\centering
	\graphicspath{{fig/}}
	\setlength{\abovecaptionskip}{0.1 pt plus 0pt minus 0pt}	   
	\includegraphics[width=0.45\textwidth]{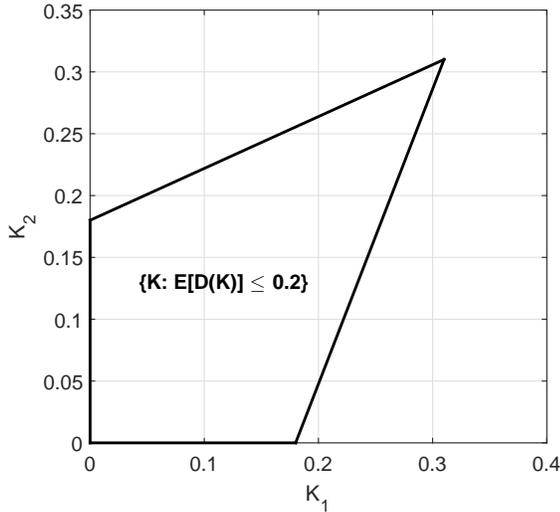}
	\figcaption{Example of Expected Maximum Drawdown Constraint}
	\label{fig: Convex_EDD}
\end{figure}

Further to the discussion of drawdown above, if the allocation vector $K$ is multi-dimensional, it would be desirable to have a convex drawdown constraint so that the log-growth optimization problem can be treated as concave program and can be solved in a very efficient way.
To this end, in this section we provide two conjectures involving convexity of the maximum drawdown.
\\ \\
For motivation, consider the single coin flipping example, it is clear from the monotonicity in Figure~\ref{fig5:E_DD_vs_K_for_motivation_ex} that expected drawdown is increasing function in~$K$. Thus, for $0<\varepsilon <1$, the $\mathbb{E}[D(K)] \leq \varepsilon$ leads to an interval constraint for $K$ which is convex.
For two identical coins for which $X_i=1$ with probability~\mbox{$p = 0.9$} and $X_i=-1$ with probability~$1-p$ for~$i=1,2$, a Monte~Carlo simulation indicates that the constraint set of expected maximum drawdown defines a convex set;~e.g., see Figure~\ref{fig: Convex_EDD}. This leads to the following conjecture.

\begin{figure}[htbp]
 	\centering
 	\graphicspath{{fig/}}
 	\setlength{\abovecaptionskip}{0.1 pt plus 0pt minus 0pt}	   
 	\includegraphics[width=0.45\textwidth]{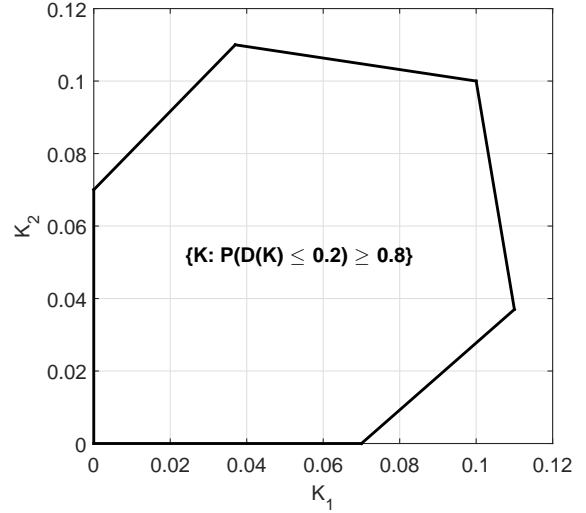}
 	\figcaption{Example of Probability of Maximum Drawdown Constraint}
 	\label{fig:Convex_prob_DD}
 \end{figure}

{\bf Conjecture 1:} {\it Given $0 < \varepsilon  < 1$, the set for the expected maximum drawdown 
$$
\{ K \in \mathcal{K}: \mathbb{E}[D(K)] \le \varepsilon \}
$$
is convex.
}
\\
\\
 Instead of constraining expected maximum drawdown, we might consider the set for the probability that the maximum drawdown stays below some level $ \varepsilon >0$. For  the same two identical coins flipping example, Figure~\ref{fig:Convex_prob_DD} shows that the drawdown constraint set is still convex. This leads to the following conjecture.\\ \\ \\
 
{\bf Conjecture 2:} {\it Given $0 < \varepsilon  < 1$ and $0< \delta < 1$, the set for the probability of maximum drawdown
$$
\{ K \in \mathcal{K}: P(D(K) \le \varepsilon ) \ge 1-\delta \}
$$is convex.
}
\\
\\
Although we have carried out several Monte Carlo simulations supporting these conjectures, in view of the fact that involving two coins, the question on the convexity is still the conjectures which may not be true for general case, we now introduce a  surrogate for expected drawdown.\\

\subsection{Surrogate for Expected Drawdown}
\begin{figure}[htbp]
	\centering
	\graphicspath{{fig/}}
	\setlength{\abovecaptionskip}{0.1 pt plus 0pt minus 0pt}	   
	\includegraphics[width=0.45\textwidth]{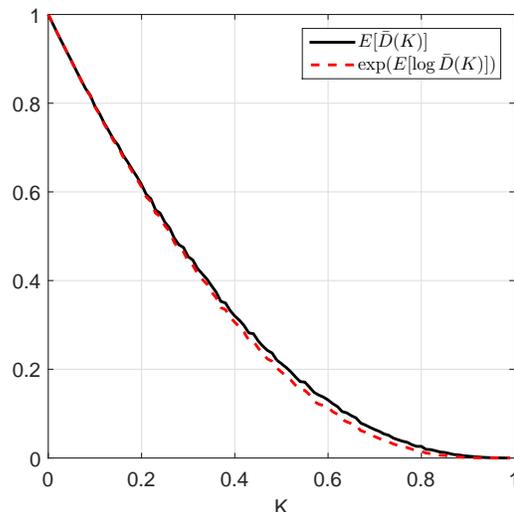}
	\figcaption{Expected Complementary Drawdown and Its Surrogate}
	\label{fig:ED_C vs K}
\end{figure}

First noting that 
\begin{align*}
	D(K) = \max_{0 \le l \le k \le N}\frac{V(l) - V(k)}{V(l)} = 1 - \bar D(K)
\end{align*}
where
$$ 
\bar D(K) \doteq \min_{0 \le l \le k \le N }\frac{ V(k)}{V(l)},
$$ define {\it complementary maximum drawdown}. We consider this complementary drawdown as a surrogate and work with the surrogate constraint 
	\[ 
	{\log  \bar D(K) } \ge \log(1- \varepsilon)
	\] 
Thus, the following lemma indicates that this complementary drawdown constraint defines a convex set.
\\
\newtheorem{lem}{\bf Lemma}[section]
{\begin{lem}\label{lemlimit} 
		{\it Given $0 <\epsilon < 1$, the set $$
			\left\{ \; K \in \mathcal{K}: \mathbb{E}[\log  \bar D(K)] \ge \log (1 - \varepsilon ) \;
			\right\}
			$$
			 is convex.
		}\\
\end{lem}}
{\bf Proof:} Given $0<\varepsilon <1$, we have
	\begin{align*}
	 \mathbb{E} [\log \bar D(K)]
	 & = \mathbb{E} \left[\log \left( {\mathop {\min }\limits_{ 0 \le l \le k \le N} \frac{{V\left( k \right)}}{{V\left( l \right)}}} \right)\right]\\
	& = \mathbb{E} \left[
		\mathop {\min }\limits_{0 \le l \le k \le N} \log \frac{{V\left( k \right)}}{{V\left( l \right)}}
	\right] \\
	& = \mathbb{E} \left[ 
	\min_{0 \le l \le k \le N} \sum\limits_{i = l}^{k - 1} \log{\left( {1 + {K^T}X\left( i \right)} \right)} 
	\right] \\
	&= \int_{\mathcal{X}}  \;\; \min_{0 \le l \le k \le N} \sum\limits_{i = l}^{k - 1} \log{\left( {1 + {K^T}x} \right)} f_X(x)dx.
	\end{align*}	
Note that the function $\sum\limits_{i = l}^{k - 1} {\log \left( {1 + {K^T}x} \right)} $ is concave in~$K$, using the fact that the minimum over an index collection of the concave functions is concave, it follows that~$\mathbb{E}[\log  \bar D(K)]$ is a concave function. Hence, the set 
	$$
	\{K \in \mathcal{K}: \mathbb{E}[\log  \bar D(K))] \ge \log(1 - \varepsilon)\}
	$$
is convex. $\square$
\\
\\
{\bf Remark:} Since log function is concave, using Jensen's inequality, we obtain
	\begin{align*}
	 \mathbb{E}[ \log \bar D(K)]& \le \log \mathbb{E}[  \bar D(K)].
	\end{align*}
Now exponentiating on both sides, we obtain
$$
	  \mathbb{E}[  \bar D(K)] \ge \exp( \mathbb{E}[ \log \bar D(K)]) .
$$
\\
To consider the tightness of this bound, we revisit the single coin flipping gamble again with probability~$p=0.9$ and~ $N=252$.  Figure~\ref{fig:ED_C vs K} provides a comparison between~$\mathbb{E}[\bar D(K)]$ and~$\exp(\mathbb{E}[ \log \bar D(K)])$ obatined by using Monte Carlo simulation. 
For this simple case, It is clear that $\mathbb{E}[\bar D(K)]$ is very close to $\exp(\mathbb{E}[ \log \bar D(K)])$. In other words, the surrogate complementary drawdown can be a drawdown candidate. \\

\section{Conclusion} 
\begin{figure}[htbp]
	\centering
	\graphicspath{{fig/}}	\setlength{\abovecaptionskip}{0.1 pt plus 0pt minus 0pt}	   
	\includegraphics[width=0.48\textwidth]{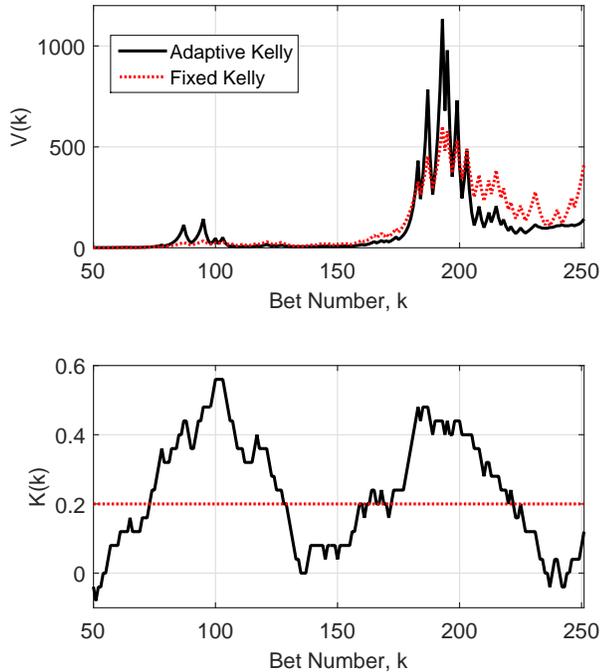}
	\figcaption{Adaptive Kelly Strategy for $M=50$ Training Bets \\}
	\label{fig:Adaptive Kelly}
\end{figure}

\parindent = 0pt
In this paper, the focal point was some of the limitations associated with application the Kelly Criterion. By way of further research, in addition to the drawdown issues described in Section~3, another possibility involves modification of the feedback control scheme~$I_i(k) = K_i V(k)$ defining the investment. Perhaps use of other variables in the ``controller'' such as the drawdown itself would result in improved performance. More generally, it would be of interest to pursue the Kelly-based theory with other risk metrics such as the Sharpe Ratio, see~\cite{Sharpe_1994}, in play.\\
\\
An important line of future research involves extension of existing results to problems involving with~$f_X(x)$ not assumed to be known. 
For example, when the theory is applied in a stock-trading context instead of assuming~$f_X(x)$ is known, it would make sense to consider the use of an adaptive scheme to obtain a~$K$-vector which is time-varying; i.e., as nature of the market dynamics change, the investment function is correspondingly adjusted.\\
\\
To provide a simple illustration how such an adaptive scheme might work, we consider the coin-flipping game described in Section~1 with initial account value $V(0) = 1$ and unknown underlying probability~$p = 0.6$. Now, the bettor, not knowing~$p$ observes outcomes~$X(k)$ and constructs a relative frequency estimate~$\hat{p}(k)$ of~$p$ using a sliding window of size~$M < N$. The first~$M$ steps constitute the training period within which no betting is done, and then, for~$k \geq M$, the estimator is given by
$$
\hat{p}(k) \doteq \frac{1}{M}\sum_{i = k-M}^{k-1} \max\{ sign(X(k)),0\}.
$$
Note that the estimator above is used for expressing that the number of winning bet.
Now, using the estimator, we can obtain investment fraction
$$
\hat{K}(k) = 2\hat{p}(k) - 1
$$
The results, summarized in Figure~\ref{fig:Adaptive Kelly}.\\

\end{document}